\documentclass[12pt]{amsart}
\usepackage{amssymb}
\usepackage{xcolor}
\usepackage{graphicx}
\usepackage{float}
\usepackage{dsfont}
\usepackage{url}
\usepackage{hyperref}
\usepackage{mathtools}
\usepackage{anysize}
\usepackage{geometry}
\marginsize{2cm}{2cm}{2cm}{2cm}
\usepackage[normalem]{ulem}
\usepackage[mathlines]{lineno}
\theoremstyle{plain}
\numberwithin{equation}{section}
\newtheorem{thm}{Theorem}[section]

\newcommand{\PP}[1]{\left(#1\right)}
\newcommand{\CC}[1]{\left[#1\right]}
\newcommand{\LL}[1]{\left\{#1\right\}}

\newcommand{\abs}[1]{\left|#1\right|}
\newcommand{\norm}[1]{\left|\left|#1\right|\right|}
\renewcommand{\Pr}[1]{\mathbb{P}\left(#1\right)}

\newcommand{\R}{\mathbb{R}}

\newcommand{\dd}{\textnormal{d}}

\newcommand{\g}{\rho_\Omega}
\newcommand{\atandos}{\operatorname{atan2}}

\title[Distribution of the angle between random segments]{Angle Distributions for Intersecting Random Segments in Star-Shaped Planar Domains}
\author{Paulo Manrique-Mir\'on}
\email{manriquemiron@gmail.com}
\subjclass[2020]{Primary 60D05; Secondary 52A30, 52A22, 53C65.}
\keywords{geometric probability, random segments, crossing angle, star-shaped domains, radial function, integral geometry, Sylvester problem}
\begin{document}
\begin{abstract}
Let $\Omega\subset\R^2$ be a bounded planar set that is star-shaped with respect to the origin, and let $A,B,C,D$ be independent random points uniformly distributed on $\Omega$. We consider the random segments $S_{AB}$ and $S_{CD}$ and study the distribution of the smaller angle $\Theta\in[0,\pi/2]$ formed by them, conditional on the event that they intersect. Using the radial function of $\Omega$, together with a parametrization of each segment in terms of its supporting line and the positions of its endpoints along that line, we derive an integral representation for the conditional distribution
\[
\Pr{\Theta\leq\theta \middle| S_{AB}\cap S_{CD}\neq\varnothing}.
\]
The resulting expression makes explicit how the geometry of the boundary of $\Omega$ determines the angular distribution. The probability of intersection appears naturally as the normalizing constant and is related to the probabilistic version of Sylvester's four-point problem.
\end{abstract}
\maketitle

\section{\textbf{Introduction}}\label{sec:intro}

Geometric probability is concerned with describing the random behavior of geometric objects generated by random variables. A classical example is Buffon's needle problem, which asks for the probability that a \emph{random} line segment intersects one of a family of equally spaced parallel lines. The solution to this problem provides a probabilistic procedure for computing the value of $\pi$, which may appear somewhat unexpected at first sight, since it amounts to counting how often a needle dropped at random onto a surface intersects a collection of equally spaced parallel lines. A wide variety of problems of this nature arise in geometric probability; see, for instance,~\cite{kendall1963geometrical,mathai1999introduction}.

In~\cite{manrique2023angle}, the problem of determining the distribution of the angle formed by two random segments whose endpoints are independent and uniformly distributed on the unit disk is studied. Despite the apparent simplicity of the underlying random mechanism, the corresponding angular distribution does not admit a simple closed-form expression. In this setting, the distribution of the angle between the two segments, conditional on their intersection, is given in Theorem~\ref{thm:mar29ago1553}.

\begin{thm}\label{thm:mar29ago1553}
For $0<\theta<\pi$, define
\begin{align}
 g^*(\theta)
 :=\frac{1}{2\pi}\sum_{\varepsilon\in\{-1,1\}}
 \int_0^1\!\int_0^1
 &g^*_{1,\varepsilon}(\rho_{AB},\rho_{CD},\theta)\nonumber\\
 &\times
 \mathds{1}_{\left\{
 \sqrt{1-\rho_{AB}^2}\,|\sin\theta|
 \geq |\rho_{AB}\cos\theta+\varepsilon\rho_{CD}|
 \right\}}
 \,\dd\rho_{AB}\,\dd\rho_{CD},
 \label{eq:corrected-gstar}
\end{align}
where
\begin{align*}
 g^*_{1,\varepsilon}(\rho_{AB},\rho_{CD},\theta)
 &:=\left(\frac{4}{\pi}\right)^2
 \sqrt{1-\rho_{AB}^2}\sqrt{1-\rho_{CD}^2}\\
 &\quad\times
 \left[
 1-\frac{\rho_{AB}^2+\rho_{CD}^2
 +2\varepsilon\rho_{AB}\rho_{CD}\cos\theta}
 {\sin^2\theta}
 \right]^2
 \mathds{1}_{\{\rho_{AB},\rho_{CD}\in[0,1]\}}.
\end{align*}
Set $g^*(0)=g^*(\pi)=0$ and
\[
 c:=\int_0^\pi g^*(s)\,\dd s.
\]
If $S_{AB},S_{CD}$ are the random segments which extremes are independent and identically distributed random variables with uniform distribution on the unit disc, then
\[
 c=\Pr{S_{AB}\cap S_{CD}\neq\varnothing}
 =\frac13\left(1-\frac{35}{12\pi^2}\right)
 =0.234826627014\ldots,
\]
and the conditional density of the small angle between $S_{AB},S_{CD}$, $\Theta$, given intersection is
\[
 g(\theta)=\frac{g^*(\theta)}{c}\,
 \mathds{1}_{\{\theta\in(0,\pi)\}}.
\]
Equivalently,
\[
\Pr{\Theta\leq\theta\mid S_{AB}\cap S_{CD}\neq\varnothing}
 =\int_0^\theta g(s)\,\dd s,
 \qquad 0\leq\theta\leq\pi.
\]
\end{thm}


It follows from Theorem~\ref{thm:mar29ago1553} that the probability that two random segments, whose endpoints are independent random points uniformly distributed on the unit disk, intersect can be explicitly determined. This setting is closely related to the classical Sylvester problem in a probabilistic framework~\cite{kabluchko2026refinement}, which concerns the probability that $n$ random points are in \emph{convex position}, that is, that their convex hull is a polytope with exactly $n$ vertices. Theorem~\ref{thm:mar29ago1553} considers the particular case of four independent random points uniformly distributed on the unit disk. In this setting, the Sylvester problem determines the normalizing constant for the density of the angle formed by the two random segments, conditional on their intersection.



Motivated by Sylvester's problem, we consider the setting in which the endpoints of the segments are independently generated random points uniformly distributed over a bounded star-shaped set $\Omega$. In \cite{manrique2023angle}, a suitable parametrization of the random segments is introduced, which makes it possible to derive expressions for the density of the angle formed by the segments, conditional on their intersection. In \cite{manrique2023angle}, the particular case in which $\Omega$ is the unit disk $\mathds{D}$ is considered.

In the present work, we introduce the {\it radial function} of $\Omega$ (see \eqref{eqn:21sep}), which encodes the geometry of the boundary of $\Omega$ and, in turn, determines the behavior of the distribution of the angle formed by two random segments conditional on their intersection.


More precisely, we assume that $\Omega\subset\mathbb{R}^2$ satisfies the following conditions: 
$0\in\Omega$, there exists $M>0$ such that $\Omega\subset B(0,M)$, where $B(0,M)$ denotes the ball centered at the origin with radius $M$, and, for every $z\in\Omega$, one has
\[
tz\in\Omega
\qquad\text{for all } t\in[0,1].
\]
In other words, the line segment joining the origin to any point $z\in\Omega$ is entirely contained in $\Omega$. Clearly, this condition is automatically satisfied whenever $\Omega$ is convex. For a set $\Omega$ satisfying these assumptions, we define its radial function $\g$ of $\Omega$ by
\begin{equation}\label{eqn:21sep}
\g(\theta):=\sup\LL{s\geq 0: s(\cos(\theta),\sin(\theta))^T\in\Omega},
\quad \theta\in[0,2\pi],
\end{equation}
where $v^T$ denotes the transpose of the vector $v$. The set $\Omega$ can be described in terms of its radial function as
\[
\Omega = \LL{(s,\theta) : 0 \leq s \leq \g(\theta)}.
\]

To illustrate the usefulness of the function $\g$, let $B(r)$ denote the ball of radius $r\geq 0$, and set
$I:=\mbox{Area}(\Omega)$. Then
\begin{align}
\frac{1}{I} \int_{\Omega} \mathds{1}_{B(r)}(w)\,\dd w
&=
\frac{1}{I}\int_{0}^{2\pi}
\int_{0}^{\g(\theta)}
\mathds{1}_{[0,r]}(s)\,s\,\dd s\,\dd\theta
\label{eqn:21sep850}\\
&=
\frac{1}{2I}
\int_{0}^{2\pi}
\min\LL{r^2,\g^2(\theta)}\,\dd\theta.
\nonumber
\end{align}

Let $X$ be a random point uniformly distributed on $\Omega$. The distribution of $\norm{X}$ can be obtained by the standard method of evaluating probabilities over small neighborhoods. Let $r>0$ and let $\Delta>0$ be sufficiently small. Then
\begin{align*}
\Pr{\norm{X}\in [r,r+\Delta]} & = \frac{1}{I} \int_{\Omega} \mathds{1}_{B(r+\Delta)}(w) dw - \frac{1}{I} \int_{\Omega} \mathds{1}_{B(r)}(w) dw \\
& = \frac{1}{2I} \int_{0}^{2\pi} \min\LL{(r+\Delta)^2,\g^2(\theta)} \dd\theta - \frac{1}{2I} \int_{0}^{2\pi} \min\LL{r^2,\g^2(\theta)} \dd\theta \\
& = \frac{1}{2I} \int_{0}^{2\pi} \mathds{1}_{\LL{(r+\Delta)^2 \geq \g^2 \geq r^2}} \CC{\g^2-r^2} \dd \theta\; + \\
& \quad\quad \frac{1}{2I} \int_{0}^{2\pi} \mathds{1}_{\LL{R^2 \geq \g^2 \geq (r+\Delta)^2}} \CC{\g^2-r^2} \dd \theta,
\end{align*} 
where, in the last line, we have set $R:=\sup_{x\in\Omega}{\norm{x}}$. Then,
\[
\lim_{\Delta \to 0} \frac{\Pr{\norm{X}\in [r,r+\Delta]}}{\Delta} = \frac{1}{I} r \int_{0}^{2\pi} \mathds{1}_{\LL{R\geq \g \geq r}} \dd\theta,
\]
Thus, the density $f_{\norm{X}}(r)$ of $\norm{X}$ is given by
\begin{equation}\label{eqn:23sep}
f_{\norm{X}}(r)
=
\frac{1}{I} r \lambda\PP{\LL{\theta\in[0,2\pi]: R\geq \g(\theta) \geq r}},
\end{equation}
where $\lambda$ denotes the Lebesgue measure.

If $\Omega$ is the unit disk, then $\g(\theta)=1$ for every 
$\theta\in[0,2\pi]$. Hence, $I=\pi$, $R=1$, and
\begin{equation}\label{eqn:23sep}
f_{\norm{X}}(r)
=
2r\,\mathds{1}_{\LL{r\in[0,1]}}.
\end{equation}

\begin{figure}
\centering
\includegraphics{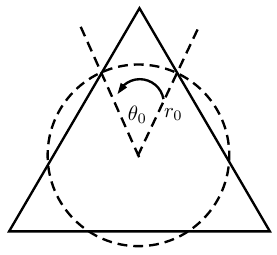}
\caption{When $\Omega$ is an equilateral triangle whose centroid is at the origin.}\label{fig01sep}
\end{figure}

Now we suppose that $\Omega$ is an equilateral triangle whose centroid is at the origin and whose side length is equal to $1$. In this case, determining the distribution of $\norm{X}$ is no longer straightforward. However, by using the definition of the radial function $\g$, one can observe that, in order to determine the quantity
\[
\lambda\PP{\LL{
\theta\in[0,2\pi]:
\frac{\sqrt{3}}{3}\geq \g(\theta)\geq r
}},
\]
it is enough to understand the relationship between $\theta_0$ and $r_0$, as illustrated in Figure~\ref{fig01sep}. This relationship is given by
\[
\theta_0
=
2\arcsin\PP{\frac{\sqrt{3}}{6r_0}}
-
\frac{\pi}{3}.
\]

The density of $\norm{X}$ in the case where $\Omega$ is an equilateral triangle with centroid at the origin is given by
\begin{equation}\label{eqn:24sep}
g_{\norm{X}}(t) = \left\{
\begin{array}{ll}
\frac{4\cdot 2\pi}{\sqrt{3}}r & \frac{\sqrt{3}}{6} \geq r \geq 0, \\
\frac{3\cdot 4}{\sqrt{3}} r\times\PP{2\arcsin\PP{\frac{\sqrt{3}}{6r}} - \frac{\pi}{3}} &\frac{\sqrt{3}}{3} \geq r >\frac{\sqrt{3}}{6},\\
0 & \mbox{otherwise}.
\end{array}
\right.
\end{equation}

Equation \eqref{eqn:23sep} provides a compact and general expression for determining the distribution of $\norm{X}$, which depends only on the geometry of the boundary of $\Omega$. As can be seen from \eqref{eqn:23sep} and \eqref{eqn:24sep}, even when the geometry of $\Omega$ appears to be simple, the behavior of the associated random object may be nontrivial.

The rest of this paper is organized as follows. In Section 2, we state the main result of this paper, which provides an expression for the density of the angle formed by two random segments whose endpoints are independent random points uniformly distributed on $\Omega$. Section 3 is devoted to the proof of the main result.

{\bf Acknowledgments}. Before proceeding, I would like to thank Michelle Zela, Marceline X., Kusy X., Yanay X., and Sophie X., for their support during the preparation of this work.

\section{Angle Between Random Segments}\label{sec:angulo}

We consider four random points $A,B,C,D$, independently generated and uniformly distributed over $\Omega$ (with the properties described in Section~\ref{sec:intro}). A random segment $S_{AB}$ is defined by
\[
S_{AB} := \LL{w\in\R^2 : w=(1-\alpha)A+\alpha B, \alpha\in[0,1]}.
\]
Let $\Theta\in[0,\pi/2]$ denote the smaller angle between $S_{AB}$ and $S_{CD}$ whenever the two segments intersect.

The quantity of interest is
\[
\Pr{\Theta \leq \theta | S_{AB}\cap S_{CD}\neq\varnothing}.
\]

To establish the main result, we define
\begin{equation}\label{eqn:25sep808}
S_{\Omega}(r,\gamma) := \LL{t:r^2 + t^2\leq \g^2(\gamma + \atandos(t,r))},
\end{equation}

\begin{equation}\label{eqn:25sep812}
K_{\Omega}(r,\gamma;q) := \int_{S_{\Omega}(r,\gamma)} \int_{S_{\Omega}(r,\gamma)} \abs{u-v} \mathds{1}_{\LL{(u-q)(v-q)\leq 0}}  \dd u \dd v,
\end{equation}
y
\begin{equation}\label{eqn:soporte}
h_{\Omega}(\phi) x:= \sup_{x\in\Omega} \langle x, u(\phi) \rangle = \max_{\theta\in[0,2\pi]} \g(\theta) \cos(\theta-\phi). 
\end{equation}

\begin{thm}\label{thm:main.thm}
If $A,B,C,D$ are independent random points uniformly distributed on $\Omega$, then
\begin{align*}
& \Pr{\Theta\leq \theta | S_{AB} \cap S_{CD} \neq \varnothing}  = \\
& \quad \frac{1}{\Pr{S_{AB} \cap S_{CD} \neq \varnothing}\times \textnormal{Area}(\Omega)^4 } \int_{0}^{2\pi}  \int_{0}^{2\pi} \mathds{1}_{\LL{\abs{\cos(\gamma_{AB}-\gamma_{CD})} \leq \sin \theta}} \\
& \quad\quad \times \int_{0}^{h_{\Omega}(\gamma_{AB})} \int_{0}^{h_{\Omega}(\gamma_{CD})} K_{\Omega}\PP{r_{AB},\gamma_{AB}; \frac{r_{AB}\cos(\gamma_{AB}-\gamma_{CD}) - r_{CD}}{\sin(\gamma_{AB}-\gamma_{CD})}} \\
& \quad\quad\quad \times K_{\Omega}\PP{r_{CD},\gamma_{CD}; \frac{r_{AB} - r_{CD}\cos(\gamma_{AB}-\gamma_{CD})}{\sin(\gamma_{AB}-\gamma_{CD})}} \dd r_{CD} \dd r_{AB} \dd\gamma_{CD}\dd\gamma_{AB}.
\end{align*}
\end{thm}

Note that the probability $\Pr{S_{AB} \cap S_{CD} \neq \varnothing}$ can be computed by observing that, for $S_{AB}$ and $S_{CD}$ to intersect, the points $A,B,C,D$ must be in convex position; that is, their convex hull must be a quadrilateral. By symmetry of the labels, it follows that
\[
\Pr{S_{AB} \cap S_{CD} \neq \varnothing} = \frac{1}{3} \Pr{\mbox{$A,B,C,D$ están en posición convexa}}.
\]

On the other hand, the points $A,B,C,D$ fail to be in convex position precisely when one of them lies in the triangle $\triangle(A,B,C)$  determined by the other three. By symmetry,
\[
\Pr{S_{AB} \cap S_{CD} \neq \varnothing} = \frac{1}{3}\CC{1-4\Pr{D\in\triangle(A,B,C)}}.
\] 


\section{Proof} \label{sec:demostracion}

Given a random point \(X\) uniformly distributed over \(\Omega\), we observe that
\[
X \stackrel{D}{=} \sqrt{R} \g(\Gamma)\PP{\begin{array}{c} \cos(\Gamma) \\ \sin(\Gamma) \end{array}},
\]
where $\stackrel{D}{=}$ means to be equal in distribution, and $R$ and $\Gamma$ are independent random variables such that $R$ is uniformly distributed on $[0,1]$ and $\Gamma$ has density
\[
\frac{\g^2(\gamma)}{\int_{0}^{2\pi} \g^2(s) \dd s} \mathds{1}_{\LL{\gamma\in[0,2\pi]}}.
\]
To verify this, it suffices to consider the Jacobian $\norm{J}$ of the transformation $(x_1,x_2)\to(r,\gamma)$ defined implicitly by
\begin{align*}
x_1 & = \sqrt{r}\g(\gamma) \cos{\gamma}, \\
x_2 & = \sqrt{r}\g(\gamma) \sin{\gamma},
\end{align*}
and recalling that $I=\mbox{Area}(\Omega)= \frac{1}{2} \int_{0}^{2\pi} \g^2(s) \dd s$, we obtain
\[
\norm{J} = \frac{1}{4\pi} \int_{0}^{2\pi} \g^2(s) \dd s = \frac{I}{2\pi}.
\]

Now, we introduce the following vectors:
\[
n(\gamma):= \PP{\cos(\gamma),\sin(\gamma)}^T,\quad u(\gamma):= \PP{-\sin(\gamma),\cos(\gamma)}^T.
\]
Consider the perpendicular from the origin $O$ to the supporting line of the segment $S_{AB}$, and denote its foot by $F$. Let $\Gamma_{AB}$ be the angle formed by $OF$ and the $x$-axis, and let $R_{AB}$ denote the perpendicular distance from $O$ to the supporting line. Furthermore, let $T_A$ and $T_B$ denote the signed distances from $F$ to the endpoints $A$ and $B$, respectively, measured along the supporting line. Then, $A,B$ can be written as
\begin{align*}
A & = R_{AB}n(\Gamma_{AB}) + T_A u(\Gamma_{AB}), \\
B & = R_{AB}n(\Gamma_{AB}) + T_B u(\Gamma_{AB}), 
\end{align*}
and analogously for the points $C,D$. This {\it reparametrization} provides a more convenient description of the segment $S_{AB}$ and, in turn, allows for a clearer analysis of the behavior of the angle between the segments $S_{AB}$ and $S_{CD}$. Once again, computing the Jacobian $\norm{J_1}$ associated with the change of variables from the joint distribution of $(R_A,\Gamma_A,R_B,\Gamma_B)$ to that of $(R_{AB},\Gamma_{AB},T_A,T_B)$, we obtain
\[
\norm{J_1}:= \abs{\frac{\partial(r_A,\gamma_A,r_B,\gamma_B)}{\partial(r_{AB},\gamma_{AB},t_A,t_B)}} = \frac{4\abs{t_A-t_B}}{\g^2(\gamma_A) \g^2(\gamma_B)}.
\]

\begin{figure}
\centering
\includegraphics{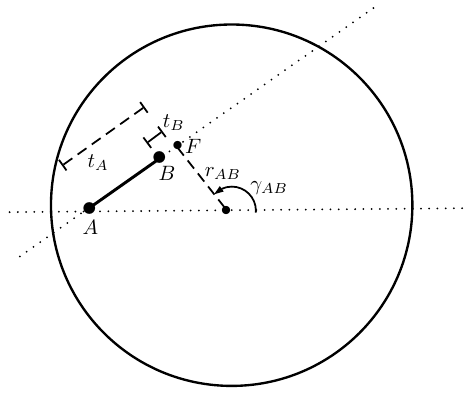}
\caption{Radial change of variables.}\label{fig:figura2}
\end{figure}

Recalling the definition of $h_{\Omega}(\phi)$,
\begin{align*}
h_{\Omega}(\phi) & := \sup_{x\in\Omega} \langle x, u(\phi) \rangle = \max_{\theta\in[0,2\pi]} \g(\theta) \cos(\theta-\phi),
\end{align*}
from the Figure \ref{fig:figura2}, we observe that
\[
0\leq r_{AB} \leq h_{\Omega}(\gamma_{AB}).
\]
We define $\operatorname{atan2}(a,b)\in(-\pi,\pi]$ as
\[
\operatorname{atan2}(a,b)
=
\begin{cases}
\arctan\left(\dfrac{a}{b}\right),
& b>0,\\[4mm]

\arctan\left(\dfrac{a}{b}\right)+\pi,
& b<0,\ a\geq 0,\\[4mm]

\arctan\left(\dfrac{a}{b}\right)-\pi,
& b<0,\ a<0,\\[4mm]

\dfrac{\pi}{2},
& b=0,\ a>0,\\[4mm]

-\dfrac{\pi}{2},
& b=0,\ a<0.
\end{cases}
\]
Geometrically, $\operatorname{atan2}(a,b)$ is the angle $\theta$ such that
\[
\PP{\cos\theta,\sin\theta}^T = \frac{1}{\sqrt{a^2+b^2}} (b,a)^T.
\]

The value $t_A$ should be satisfied that
\[
r_{AB}^2 + t_A^2 \leq \g^2\PP{\gamma_{AB} + \atandos{t_A,r_{AB}}}.
\]

It follows from the above that the joint density of
$(R_{AB},\Gamma_{AB},T_A,T_B)$ is given by
\begin{align}
f_{(R_{AB},\Gamma_{AB},T_A,T_B)}
(r_{AB},\gamma_{AB},t_A,t_B)
&=
\frac{\g^2(\gamma_A)\g^2(\gamma_B)}{4I^2}
\times
\frac{4\abs{t_A-t_B}}
{\g^2(\gamma_A)\g^2(\gamma_B)}
\nonumber\\
&=
\frac{1}{\operatorname{Area}(\Omega)^2}
\abs{t_A-t_B}.
\label{eqn:densidad}
\end{align}
on the domain
\[
\begin{array}{ll}
\mathcal{D}_{\Omega} := \left\{ (r_{AB},\gamma_{AB},t_A,t_B) \right. : & 0\leq \gamma_{AB} \leq 2\pi, \\ [1ex]
& 0\leq r_{AB} \leq h_{\Omega}(\gamma_{AB}), \\ [1ex]
& r_{AB}^2 + t_A^2 \leq \g^2\PP{\gamma_{AB} + \atandos\PP{t_A,r_{AB}}}, \\ [1ex]
& \left. r_{AB}^2 + t_B^2 \leq \g^2\PP{\gamma_{AB} + \atandos\PP{t_B,r_{AB}}} \right\}.
\end{array}
\]
In \cite{manrique2023angle}, the quantity $\delta:=\gamma_{AB}-\gamma_{CD}$ is introduced. Since $A,B,C,D$ are continuously distributed random variables, $\delta\neq 0$ almost surely. The following quantities are also defined:
\begin{align*}
q_{AB} & := \frac{R_{AB}\cos\delta - R_{CD}}{\sin \delta},\\
q_{CD} & := \frac{R_{AB} - R_{CD}\cos\delta}{\sin \delta}.
\end{align*}
Now, we define the events
\[
\mathcal{I}_{AB}:=\LL{(T_A - q_{AB})(T_B-q_{AB}) \leq 0},\quad \mathcal{I}_{CD}:=\LL{(T_C - q_{CD})(T_D - q_{CD})\leq 0},
\]
hence
\[
\mathcal{I}_{AB} \cap \mathcal{I}_{CD} \iff S_{AB}\cap S_{CD} \neq\varnothing.
\]
If $\Theta$ denotes the smaller angle between $S_{AB}$ and $S_{CD}$ whenever the two segments intersect, then
\[
\Theta = \arccos\PP{\abs{\cos\delta}} \in [0,\pi/2],
\]
i.e.,
\[
\LL{\Theta\leq \theta}=\LL{\abs{\cos\delta}\leq\sin \theta}.
\]
Therefore, to determine the density under consideration, it suffices to analyze the quantity
\[
\Pr{\abs{\cos\delta}\leq\sin\theta, \mathcal{I}_{AB},\mathcal{I}_{CD} }.
\]
From the definition of
\begin{equation*}
S_{\Omega}(r,\gamma) := \LL{t:r^2 + t^2\leq \g^2(\gamma + \atandos(t,r))},
\end{equation*}
and
\begin{equation*}
K_{\Omega}(r,\gamma;q) := \int_{S_{\Omega}(r,\gamma)} \int_{S_{\Omega}(r,\gamma)} \abs{u-v} \mathds{1}_{\LL{(u-q)(v-q)\leq 0}}  \dd u \dd v,
\end{equation*} 
we have that
\begin{align*}
& \Pr{\cos\delta\leq\sin \theta, \mathcal{I}_{AB},\mathcal{I}_{CD} }  = \\
& \quad \frac{1}{\mbox{Area}(\Omega)^4} \int_{0}^{2\pi}  \int_{0}^{2\pi} \mathds{1}_{\LL{\abs{\cos(\gamma_{AB}-\gamma_{CD})} \leq \sin \theta}} \\
& \quad\quad \times \int_{0}^{h_{\Omega}(\gamma_{AB})} \int_{0}^{h_{\Omega}(\gamma_{CD})} K_{\Omega}\PP{r_{AB},\gamma_{AB}; \frac{r_{AB}\cos(\gamma_{AB}-\gamma_{CD}) - r_{CD}}{\sin(\gamma_{AB}-\gamma_{CD})}} \\
& \quad\quad\quad \times K_{\Omega}\PP{r_{CD},\gamma_{CD}; \frac{r_{AB} - r_{CD}\cos(\gamma_{AB}-\gamma_{CD})}{\sin(\gamma_{AB}-\gamma_{CD})}} \dd r_{CD} \dd r_{AB} \dd\gamma_{CD}\dd\gamma_{AB}.
\end{align*}
When $\theta=\frac{\pi}{2}$, the resulting quantity is $\Pr{\mathcal{I}_{AB},\mathcal{I}_{CD} }$, which is the normalizing constant for the conditional distribution $\Pr{\Theta \leq \theta | S_{AB}\cap S_{CD}\neq\varnothing}$. This proof Theorem~\ref{thm:main.thm}.

Note that when \(\Omega\) is convex, the set \(S_{\Omega}(r,\gamma)\), which describes the admissible values of \(t\), can be identified with the intersection of \(\Omega\) and the line
\[
L(r,\gamma) := \LL{ru(\gamma) + t v(\gamma) : t \in \R},
\]
For fixed \(r\) and \(\gamma\), the intersection \(\Omega\cap L(r,\gamma)\) is a line segment. Consequently, the admissible values of \(t\) form a closed interval, which simplifies the expression for \(K_{\Omega}(r,\gamma;q)\).

\end{document}